\documentclass[12pt, twoside, eqno]{article}
\usepackage{latexsym}
\usepackage{amssymb}
\usepackage{amsfonts}
\begin{document}

\noindent {\bf \large On intra-regular and left regular and left duo
ordered\\$\Gamma$-semigroups}\bigskip

\medskip

\noindent{\bf Niovi Kehayopulu, Michael Tsingelis}\bigskip

\bigskip{\small

\noindent{\bf Abstract.} For an intra-regular or a left regular and 
left duo ordered $\Gamma$-semigroup $M$, we describe the principal 
filter of $M$ which plays an essential role in the structure of this 
type of $po$-$\Gamma$-semigroups. We also prove that an ordered 
$\Gamma$-semigroup $M$ is intra-regular if and only if the ideals of 
$M$ are semiprime and it is left (right) regular and left (right) duo 
if and only if the left (right) ideals of $M$ are semiprime.\medskip

\medskip

}
\section{Introduction and prerequisites}Croisot, who used the term 
``inversive" instead of ``regular", connects the matter of 
decomposition of a semigroup with the regularity and semiprime 
conditions [1]. A semigroup $S$ is said to be left (resp. right) 
regular if for every $a\in S$ there exists $x\in S$ such that 
$a=xa^2$ (resp. $a=a^2x$). That is, if $a\in Sa^2$ (resp. $a\in 
a^2S)$ for every $a\in S$ which is equivalent to saying that 
$A\subseteq A^2S$ (resp. $A\subseteq SA^2$) for every $A\subseteq S$. 
A semigroup $S$ is said to be intra-regular if for every $a\in S$ 
there exist $x,y\in S$ such that $a=xa^2 y$. In other words, if $a\in 
Sa^2 S$ for every $a\in S$ or $A\subseteq SA^2S$ for every 
$A\subseteq S$. For decompositions of an intra-regular, left regular 
or both left regular and right regular semigroup we refer to [2, 10]. 
The concepts of intra-regular ordered semigroup and of right regular 
ordered semigroups have been introduced in [3, 4] in which the 
decomposition of an intra-regular ordered semigroup into simple 
components and the decomposition of a right regular and right duo 
ordered semigroup into right simple components has been studied. The 
principal filter of $S$ has a very simple form, both for ordered and  
non-ordered case of semigroups, and it plays an essential role in the 
decompositions.

For two nonempty sets $M$ and $\Gamma$, we denote by $A\Gamma B$ the 
set containing the elements of the form $a\gamma b$ where $a\in A$, 
$\gamma\in\Gamma$ and $b\in B$. That is, we define $$A\Gamma 
B:=\{a\gamma b \mid a\in A, b\in B, \gamma\in\Gamma\}.$$Then $M$ is 
called a {\it $\Gamma$-semigroup} if the following assertions are 
satisfied: \begin{enumerate}
\item[(1)] $M\Gamma M\subseteq M$;
\item[(2)] $a\gamma (b\mu c)=(a\gamma b)\mu c$ for all $a,b,c\in M$ 
and all $\gamma,\mu\in\Gamma$;
\item[(3)] if $a,b,c,d\in M$ and $\gamma,\mu\in\Gamma$ such that 
$a=c$, $\gamma=\mu$ and $b=d$, then $a\gamma b=c\mu 
d$.\end{enumerate}
An {\it ordered $\Gamma$-semigroup} (shortly, {\it 
$po$-$\Gamma$-semigroup}) is a $\Gamma$-semigroup $M$ with an order 
relation ``$\le$" on $M$ such that $a\le b$ implies $ac\le bc$ and 
$ca\le cb$ for every $c\in M$. A nonempty subset $A$ of $M$ is called 
a {\it subsemigroup} of $M$ if, for every $a,b\in A$ and every 
$\gamma\in\Gamma$, we have $a\gamma b\in A$. A subsemigroup $F$ of 
$M$ is called a {\it filter} of $M$ if
(1) for every $a,b\in F$ and every $\gamma\in\Gamma$ such that 
$a\gamma b\in F$, we have $a\in F$ and $b\in F$ and (2) if $a\in F$ 
and $M\ni b\ge a$, then $b\in F$. For an element $x$ of $M$, we 
denote by $N(x)$ {\it the filter of $M$ generated by $x$} (that is, 
the least with respect to the inclusion relation filter of $M$ 
containing $x$). A nonempty subset $A$ of $M$ is called a {\it left} 
(resp. {\it right}) {\it ideal} of $M$ if (1) $M\Gamma A\subseteq A$ 
(resp. $A\Gamma M\subseteq A$) and (2) if $a\in A$ and $M\ni b\le a$, 
then $b\in A$. It is called an {\it ideal} or ({\it two-sided ideal}) 
of $M$ if it is both a left and right ideal of $M$. A 
$po$-$\Gamma$-semigroup $M$ is called {\it left} (resp. {\it right}) 
{\it duo} if the left (resp. right) ideals of $M$ are two-sided. A 
subset $T$ of $M$ is called {\it semiprime} if for every $x\in M$ and 
every $\gamma\in\Gamma$ such that $x\gamma x\in T$, we have $x\in T$. 
For a subset $H$ of $M$ we denote by $(H]$ the subset of $M$ defined 
by $(H]=\{t\in M \mid t\le a \mbox { for some } t\in H\}$. We clearly 
have $M=(M]$, and for any subsets $A, B, C, D$ of $M$, we have 
$A\subseteq (A]=((A]]$; if $A\subseteq B$, then $(A]\subseteq (B]$; 
if $A\subseteq B$ and $C\subseteq D$, then $(A\Gamma C]\subseteq 
(B\Gamma D]$; $(A]\Gamma (B]\subseteq (A\Gamma B]$; and $((A]\Gamma 
(B]]=((A]\Gamma B]=(A\Gamma (B]]=(A\Gamma B]$. As we know, some 
results on semigroups (ordered semigroups) can be transferred into 
$\Gamma$-semigroups ($po$-$\Gamma$-semigroups) just putting a Gamma 
in the appropriate place, while for some other results the transfer 
is not easy.  A $\Gamma$-semigroup $M$ is called {\it intra-regular} 
if $a\in M\Gamma a\Gamma a\Gamma M$ for every $a\in M$, equivalently 
if $A\subseteq M\Gamma A\Gamma A\Gamma M$ for every $A\subseteq M$. 
It is called {\it left} (resp. {\it right}) {\it regular} if $a\in 
M\Gamma a\Gamma a$ (resp. $a\in \Gamma a\Gamma M$) for every $a\in M$ 
or $A\subseteq M\Gamma A\Gamma A$ (resp. $A\subseteq \Gamma A\Gamma 
M$) for every $A\subseteq M$. An ordered $\Gamma$-semigroup $M$ is 
called {\it intra-regular} if for every $a\in M$ we have $a\in 
(M\Gamma a\Gamma a\Gamma M]$, equivalently if for every $A\subseteq 
M$ we have $A\subseteq (M\Gamma A\Gamma A\Gamma M]$. An ordered 
$\Gamma$-semigroup $M$ is called {\it left} (resp. {\it right}) {\it 
regular} if $a\in (M\Gamma a\Gamma a]$ (resp. $(a\in a\Gamma a\Gamma 
M]$) for every $a\in M$, equivalently if $A\subseteq (M\Gamma A\Gamma 
A]$ (resp. $A\subseteq (A\Gamma A\Gamma M]$) for every $A\subseteq 
M$. Although some interesting results on $\Gamma$-semigroups are 
obtained with these definitions, these definitions fail to describe 
the principal filter of intra-regular, left regular and right regular 
$\Gamma$-semigroups (ordered $\Gamma$-semigroups) which play an 
essential role in the investigation. To overcome this difficulty, in 
[8] a new definition of intra-regular and of left regular 
$\Gamma$-semigroups has been introduced. The intra-regular 
$\Gamma$-semigroup has been defined as a $\Gamma$-semigroup such that 
$a\in M\Gamma a\gamma a\Gamma M$ for each $a\in M$ and each 
$\gamma\in\Gamma$ and the left (resp. right) regular 
$\Gamma$-semigroup as a $\Gamma$-semigroup in which $a\in M\Gamma 
a\gamma a$ (resp. $a\in a\gamma a\Gamma M$) for each $a\in M$ and 
each $\gamma\in\Gamma$ and it is proved that a $\Gamma$-semigroup $M$ 
is left regular (in that new sense) if and only if it is a union of a 
family of left simple subsemigroups on $M$. And in [9] we gave some 
further structure theorems of this type of $\Gamma$-semigroups using 
that new definition and the form of principal filters. But what 
happens in case of intra-regular and left or right regular 
$po$-$\Gamma$-semigroups? Can we describe the form of the principal 
filters using some new definitions like in the non-ordered case? The 
present paper gives the related answer. For more information on 
$\Gamma$ (or $po$-$\Gamma$)-semigroups cf., for example, the papers 
in [5--7] of the References, and the papers in which these papers 
refer. Examples on $\Gamma$-semigroups are also given in these 
papers.
\section{On intra-regular ordered $po$-$\Gamma$-semigroups}
We characterize here the intra-regular $po$-$\Gamma$-semigroups in 
terms of filters, and we prove that a $po$-$\Gamma$-semigroup $M$ is 
intra-regular if and only if the ideals of $M$ are 
semiprime.\medskip

\noindent{\bf Definition 1.} An ordered $\Gamma$-semigroup $M$ is 
called {\it intra-regular} if$$x\in (M\Gamma x\gamma x\Gamma M]$$for 
every $x\in M$ and every  $\gamma\in\Gamma$.\medskip

\noindent{\bf Theorem 2.} {\it An ordered $\Gamma$-semigroup $M$ is 
intra-regular if and only if, for every $x\in M$, we have $$N(x)=
\{y\in M \mid x\in (M\Gamma y\Gamma M]\}.$$}{\bf Proof.}
$\Longrightarrow$. Let $x\in M$ and $T:=
\{y\in M \mid x\in (M\Gamma y\Gamma M]\}$. Then we have the 
following:

(1) $T$ is a nonempty subset of $M$. Indeed: Take an element 
$\gamma\in\Gamma$ ($\Gamma\not=\emptyset$). Since $M$ is 
intra-regular, we have
$$x\in (M\Gamma x\gamma x\Gamma M]={\Big(}(M\Gamma x)\gamma x\Gamma 
M{\Big]}\subseteq {\Big(}(M\Gamma M)\Gamma x\Gamma M{\Big]}\subseteq 
(M\Gamma x\Gamma M],$$ so $x\in T$.

(2) Let $a,b\in T$ and $\gamma\in\Gamma$. Then $a\gamma b\in T$. 
Indeed: Since $a\in T$, we have $x\in (M\Gamma a\Gamma M]$. Since 
$b\in T$, we have $x\in (M\Gamma b\Gamma M]$. Since $M$ is 
intra-regular, $x\in M$ and $\gamma\in\Gamma$, we have $x\in (M\Gamma 
x\gamma x\Gamma M]$. Then we have\begin{eqnarray*}x\in (M\Gamma 
x\gamma x\Gamma M]&\subseteq&{\Big(}M\Gamma (M\Gamma b\Gamma M]\gamma 
(M\Gamma a\Gamma M]\Gamma M{\Big]}\\&=&{\Big(}M\Gamma (M\Gamma 
b\Gamma M)\gamma (M\Gamma a\Gamma M)\Gamma 
M{\Big]}\\&=&{\Big(}(M\Gamma M)\Gamma (b\Gamma M\gamma M\Gamma 
a)\Gamma (M\Gamma M){\Big]}\\&\subseteq&{\Big(}M\Gamma (b\Gamma 
M\gamma M\Gamma a)\Gamma M{\Big]}.\end{eqnarray*}We prove that 
$b\Gamma M\gamma M\Gamma a\subseteq {\Big(}M\Gamma (a\gamma b)\Gamma 
M{\Big]}$. Then we have\begin{eqnarray*}x&\in& {\bigg(}M\Gamma 
{\Big(}M\Gamma (a\gamma b)\Gamma M{\Big]}\Gamma 
M{\bigg]}={\bigg(}M\Gamma {\Big(}M\Gamma (a\gamma b)\Gamma 
M{\Big)}\Gamma M{\bigg]}\\&=&{\Big(}(M\Gamma M)\Gamma (a\gamma 
b)\Gamma (M\Gamma M){\Big]}\subseteq {\Big(}M\Gamma (a\gamma b)\Gamma 
M{\Big]},\end{eqnarray*} so $a\gamma b\in T$. Let now $b\lambda 
u\gamma v\delta a\in b\Gamma M\gamma M\Gamma a$ for some $u,v\in M$, 
$\lambda,\delta\in\Gamma$. Since $M$ is intra-regular, for the 
elements $b\lambda u\gamma v\delta a\in M$ and $\gamma\in\Gamma$, we 
have\begin{eqnarray*}b\lambda u\gamma v\delta a&\in& {\Big(}M\Gamma 
(b\lambda u\gamma v\delta a)\gamma (b\lambda u\gamma v\delta a)\Gamma 
M{\Big]}\\&=&{\Big(}(M\Gamma b\lambda u\gamma v)\delta (a\gamma 
b)\lambda (u\gamma v\delta a\Gamma 
M){\Big]}\\&\subseteq&{\Big(}M\Gamma (a\gamma b)\Gamma M{\Big]}.
\end{eqnarray*}

(3) Let $a,b\in M$ and $\gamma\in\Gamma$ such that $a\gamma b\in T$. 
Then $a,b\in T$. Indeed: Since $a\gamma b\in T$, we have $x\in 
{\Big(}M\Gamma (a\gamma b)\Gamma M{\Big]}\subseteq {\Big(}M\Gamma 
a\gamma (M\Gamma M){\Big]}\subseteq (M\Gamma a\Gamma M]$, so $a\in 
T$. Since $x\in {\Big(}M\Gamma (a\gamma b)\Gamma M{\Big]}\subseteq 
{\Big(}(M\Gamma M)\gamma b\Gamma M{\Big]}\subseteq (M\Gamma b\Gamma 
M]$, we have $b\in T$.

(4) Let $a\in T$ and $M\ni b\ge a$. Then $b\in T$. Indeed: Since 
$a\in T$, we have $x\in (M\Gamma a\Gamma M]$. Since $a\le b$, we have 
$(M\Gamma a\Gamma M]\subseteq (M\Gamma b\Gamma M]$. Then we have 
$x\in (M\Gamma b\Gamma M]$, and $b\in T$.

(5) Let $F$ be a filter of $M$ such that $x\in F$. Then $T\subseteq 
F$. Indeed: Let $a\in T$. Then $x\in (M\Gamma a\Gamma M]$, so $F\ni 
x\le u\lambda (a\mu v)$ for some $u,v\in M$, $\lambda,\mu\in \Gamma$. 
Since $F$ is a filter of $M$, $x\in F$ and $M\ni u\lambda (a\mu v)\ge 
x$, we have $u\lambda (a\mu v)\in F$. Since $F$ is a filter of $M$, 
$u, a\mu v\in M$, $\lambda\in \Gamma$ and $u\lambda (a\mu v)\in F$, 
we have $a\mu v\in F$, again since $F$ is a filter of $M$, $a,v\in M$ 
and $\mu\in \Gamma$, we have $a\in F$. \\
\noindent$\Longleftarrow$. Let $x\in M$ and $\gamma\in\Gamma$. Then
$x\in (M\Gamma x\gamma x\Gamma M]$. Indeed: Since $N(x)$ is a 
subsemigroup of $M$, $x\in N(x)$ and $\gamma\in\Gamma$, we have 
$x\gamma x\in N(x)$. By hypothesis, we get  
$x\in{\Big(}M\Gamma(x\gamma x)\Gamma M{\Big]}=(M\Gamma x\gamma 
x\Gamma M]$, thus $M$ is intra-regular.$\hfill\Box$\medskip

\noindent{\bf Theorem 3.} {\it An ordered $\Gamma$-semigroup $M$ is 
intra-regular if and only if the ideals of $M$ are 
semiprime}.\medskip

\noindent{\bf Proof.} $\Longrightarrow$. Let $A$ be an ideal of $M$, 
$x\in M$ and $\gamma\in\Gamma$ such that $x\gamma x\in A$. Since $M$ 
is intra-regular, we have $x\in {\Big(}M\Gamma (x\gamma x)\Gamma 
M{\Big]}\subseteq {\Big(}(M\Gamma A)\Gamma M{\Big]}\subseteq (A\Gamma 
M]\subseteq (A]=A$, then $x\in A$, and $A$ is semiprime.\\
$\Longleftarrow$. Let $x\in M$ and $\gamma\in\Gamma$. Then $x\in 
(M\Gamma x\gamma x\Gamma M]$. In fact: The set $(M\Gamma x\gamma 
x\Gamma M]$ is an ideal of $M$. This is because it is a nonempty 
subset of $M$, $M\Gamma (M\Gamma x\gamma x\Gamma M]\subseteq 
{\Big(}M\Gamma(M\Gamma x\gamma x\Gamma 
M]{\Big]}={\Big(}M\Gamma(M\Gamma x\gamma x\Gamma M){\Big]}\subseteq 
(M\Gamma x\gamma x\Gamma M]$, $(M\Gamma x\gamma x\Gamma M]\Gamma 
M\subseteq (M\Gamma x\gamma x\Gamma M]$, and ${\Big(}(M\Gamma x\gamma 
x\Gamma M]{\Big]}=(M\Gamma x\gamma x\Gamma M]$. Since $(M\Gamma 
x\gamma x\Gamma M]$ is semiprime and
$(x\gamma x)\gamma (x\gamma x)=x\gamma (x\gamma x)\gamma x\in M\Gamma 
x\gamma x\Gamma M\subseteq (M\Gamma x\gamma x\Gamma M]$, we have 
$x\gamma x\in (M\Gamma x\gamma x\Gamma M]$. Again since $(M\Gamma 
x\gamma x\Gamma M]$ is semiprime, we have $x\in (M\Gamma x\gamma 
x\Gamma M]$, so $M$ is intra-regular.$\hfill\Box$
\section{On left regular and left duo $po$-$\Gamma$-semigroups}First 
we notice that the left (and the right) $po$-$\Gamma$-semigroups are 
intra-regular. Then we characterize the $po$-$\Gamma$-semigroups 
which are both left regular and left duo in terms of filters and we 
prove that a $po$-$\Gamma$-semigroup $M$ is left (resp. right) 
regular if and only if the left (resp. right) ideals of $M$ are 
semiprime.\medskip

\noindent{\bf Definition 4.} An ordered $\Gamma$-semigroup $M$ is 
called {\it left regular} (resp. {\it right regular}) if$$x\in 
(M\Gamma x\gamma x] \mbox { (resp. } x\in (x\gamma x\Gamma M])$$for 
every $x\in M$ and every  $\gamma\in\Gamma$.\medskip

\noindent{\bf Proposition 5.} {\it Let M be an ordered 
$\Gamma$-semigroup. If M is left (resp. right) regular, then M is 
intra-regular}.\medskip

\noindent{\bf Proof.} Let $M$ be left regular, $x\in M$ and 
$\gamma\in\Gamma$. Then we have\begin{eqnarray*}x\in (M\Gamma x\gamma 
x]&\subseteq&{\Big(}M\Gamma (M\Gamma x\gamma x]\gamma x{\Big 
]}={\Big(}M\Gamma (M\Gamma x\gamma x)\gamma x{\Big ]}\\&\subseteq 
&{\Big(}(M\Gamma M)\Gamma (x\gamma x)\Gamma M{\Big]}\subseteq 
{\Big(}M\Gamma x\gamma x\Gamma M{\Big]},\end{eqnarray*}thus $M$ is 
intra-regular.$\hfill\Box$\medskip

\noindent{\bf Theorem 6.} {\it An ordered $\Gamma$-semigroup $M$ is 
left regular and left duo if and only if, for every $x\in M$, we have 
$$N(x)=
\{y\in M \mid x\in (M\Gamma y]\}.$$}{\bf Proof.}
$\Longrightarrow$. Let $x\in M$ and $T:=
\{y\in M \mid x\in (M\Gamma y]\}$. Since $M$ is left regular, we 
have
$x\in (M\Gamma x\gamma x]\subseteq {\Big(}(M\Gamma M)\Gamma 
x{\Big]}\subseteq (M\Gamma x]$, so $x\in T$, and $T$ is a nonempty 
subset of $M$. Let $a,b\in T$ and $\gamma\in\Gamma$. Since $x\in 
(M\Gamma a]$, $x\in (M\Gamma b]$ and $M$ is left regular, we 
have\begin{eqnarray*}x\in (M\Gamma x\gamma 
x]&\subseteq&{\Big(}M\Gamma (M\Gamma b]\gamma (M\Gamma 
a]{\Big]}={\Big(}M\Gamma (M\Gamma b)\gamma (M\Gamma 
a){\Big]}\\&\subseteq&{\Big(}M\Gamma (b\gamma M\Gamma 
a){\Big]}.\end{eqnarray*}In addition, $b\gamma M\Gamma a\subseteq 
(M\Gamma a\gamma b]$. Indeed: Let $b\gamma u\mu a\in b\gamma M\Gamma 
a$, where $u\in M$ and $\mu\in\Gamma$. Since $M$ is left regular, we 
have $$b\gamma u\mu a\in {\Big(}M\Gamma (b\gamma u\mu a)\gamma 
(b\gamma u\mu a){\Big]}\subseteq {\Big(}M\Gamma (a\gamma b) \Gamma 
M{\Big]}={\Big(}(M\Gamma a\gamma b]\Gamma M{\Big]}.$$Since $(M\Gamma 
a\gamma b]$ is a left ideal, it is a right ideal of $M$ as well, so 
$(M\Gamma a\gamma b]\Gamma M\subseteq (M\Gamma a\gamma b]$, then 
$b\gamma u\mu a\in{\Big(}(M\Gamma a\gamma b]{\Big]}=(M\Gamma a\gamma 
b]$. Hence we obtain$$x\in {\Big(}M\Gamma(M\Gamma a\gamma 
b]{\Big]}={\Big(}M\Gamma(M\Gamma a\gamma b){\Big]}\subseteq 
{\Big(}M\Gamma(a\gamma b){\Big]},$$from which $a\gamma b\in T$.\\Let 
$a,b\in M$ and $\gamma\in\Gamma$ such that $a\gamma b\in T$. Since 
$x\in (M\Gamma a\gamma b]\subseteq (M\Gamma b]$, we have $b\in T$. 
Besides, $x\in (M\Gamma a\gamma b]\subseteq {\Big(}(M\Gamma a]\Gamma 
M{\Big]}$. The set $(M\Gamma a]$ as a left ideal, is a right ideal of 
$M$ as well, so $(M\Gamma a]\Gamma M\subseteq (M\Gamma a]$. Thus we 
have $x\in {\Big(}(M\Gamma a]{\Big]}=(M\Gamma a]$, and $a\in T$.\\Let 
$a\in T$ and $M\ni b\ge a$. Since $M$ is left regular, we have $$x\in 
(M\Gamma a\gamma a]\subseteq (M\Gamma b\gamma b]\subseteq 
{\Big(}(M\Gamma b]\Gamma M{\Big]}.$$$(M\Gamma b]$ as a left ideal is 
a right ideal of $M$, so $(M\Gamma b]\Gamma M\subseteq (M\Gamma b]$. 
Hence we have $x\in {\Big(}(M\Gamma b]{\Big]}=(M\Gamma b]$, and $b
\in T$.\\Let $F$ be a filter of $M$ such that $x\in F$ and let $a\in 
T$. Since $x\in (M\Gamma a]$, we have $F\ni x\le u\mu a$ for some 
$u\in M$, $\mu\in\Gamma$. Since $F$ is a filter of $M$, we have $u\mu 
a\in F$, and $a\in F$.\\$\Longleftarrow$. Let $x\in M$ and 
$\gamma\in\Gamma$. Since $x\in N(x)$ and $N(x)$ is a subsemigroup of 
$M$, we have $x\gamma x\in N(x)$. By hypothesis, we get $x\in 
(M\Gamma x\gamma x]$, so $M$ is left regular. Let now $A$ be a left 
ideal of $M$, $a\in A$, $\gamma\in\Gamma$ and $u\in M$. Since 
$a\gamma u\in N(a\gamma u)$ and $N(a\gamma u)$ is a filter of $M$, we 
have $a\in N(a\gamma u)$. By hypothesis, we have $a\gamma u\in 
(M\Gamma a]\subseteq (M\Gamma A]\subseteq (A]=A$. Thus $A$ is right 
ideal of $M$. $\hfill\Box$\\The right analogue of Theorem 6 also 
holds, and we have\medskip

\noindent{\bf Theorem 7.} {\it An ordered $\Gamma$-semigroup $M$ is 
right regular and right duo if and only if, for every $x\in M$, we 
have $$N(x)=\{y\in M \mid x\in (y\Gamma M]\}.$$}
\noindent{\bf Theorem 8.} {\it An ordered $\Gamma$-semigroup $M$ is 
left (resp. right) regular if and only if the left (resp. right) 
ideals of $M$ are semiprime}.\medskip

\noindent{\bf Proof.} $\Longrightarrow$. Let $M$ be left regular, $A$ 
a left ideal of $M$, $x\in M$ and $\gamma\in\Gamma$ such that 
$x\gamma x\in A$. Then we have $x\in {\Big(}M\Gamma (x\gamma 
x){\Big]}\subseteq (M\Gamma A]\subseteq (A]=A$, so $M$ is 
semiprime.\\$\Longleftarrow$. Suppose the left ideals of $M$ are 
semiprime and let $x\in M$ and $\gamma\in\Gamma$. Since $(M\Gamma 
x\gamma x]$ is a left ideal of $M$ and $(x\gamma x)\gamma (x\gamma 
x)\in (M\Gamma x\gamma x]$, we have $x\gamma x\in (M\Gamma x\gamma 
x]$, and $x\in (M\Gamma x\gamma x]$, so $M$ is left 
regular.$\hfill\Box$\\{\small

\end{document}